\documentclass[a4paper, leqno, 12pt]{amsart}
 \usepackage[usenames]{color}
 \usepackage{amssymb, amsthm, amsmath, amscd}
 \usepackage[dvipdfmx]{graphicx}
\usepackage{tikz} 
\usetikzlibrary{positioning, intersections, calc, arrows.meta,math} 
\usepackage{ascmac}
 \usepackage{fancybox}
 \usepackage{ulem}
 \usepackage[all]{xy}
 \usepackage{comment}
 \usepackage{here}
 \usepackage{float}
 \input amssym.def
 \input amssym.tex
 \newcommand{\Z}{\mathbb{Z}}

 \newcommand{\C}{\mathbb{C}}

  \newcommand{\Ov}[1]{\overline{#1}}
 \newcommand{\Ti}[1]{\widetilde{#1}}
  
  \newcommand{\Bb}[1]{\mathbb{#1}}
  
   \newcommand{\Rm}[1]{\mathrm{#1}}
   \newcommand{\Ca}[1]{\mathcal{#1}} 
\newtheorem{prop}{Proposition}[section]
 \newtheorem{thm}[prop]{Theorem}
 \newtheorem{lem}[prop]{Lemma}

\theoremstyle{remark}
\newtheorem{define}[prop]{Definition}
\newtheorem{rmk}[prop]{Remark}
  \newtheorem{exa}[prop]{Example}
    \newtheorem{prob}[prop]{Problem}

\numberwithin{equation}{section} 
 
 \makeatletter
\@namedef{subjclassname@2020}{2020 Mathematical Subject Classification}
\makeatother

\begin{document}
\baselineskip=18pt
 \title{On the Gonality type invariants and the slope of a fibered $3$-fold}
\author{Hiroto Akaike}

 
 
\maketitle

\begin{abstract}
The \textit{slope} of a fibered $3$-folds $f:X \to B$ is a relative numerical invariant defined by
$\lambda(f) := K_{f}^{3}/\mathrm{deg}(f_{\ast}\omega_{f})$,
where $K_{f}$ is the relative canonical divisor and $\omega_{f}$ is the relative dualizing sheaf.
Establishing slope inequalities is a fundamental problem in the geography of fibered spaces. 
In this paper, we introduce a new invariant called the \textit{minimal covering degree} as a gonality-type invariant and study a lower bound of the slope increasing with the covering gonality and the minimal covering degree of the general fiber of $f$.
\end{abstract}

\section{Introduction}
We work out over the complex number field $\C$.
In this paper, 
a {\it $($algebraic$)$ fiber space} $f:X \to B$ is a surjective morphism 
from a normal projective variety $X$ to a smooth projective curve $B$
with connected fibers.
We assume that a general fiber $F$ of $f$ is of general type.
If $X$ is a smooth projective surface,
we call it a {\it fibered surface} and if $X$ is a $3$-fold, we call it a {\it fibered $3$-fold}.

\subsection{Gonality type invariants}

In order to state our main results,
we need to introduce the two gonality type invariants, {\it covering gonality} and {\it minimal covering degree}.
To define covering gonality, 
we first introduce a {\it covering family of curves} on $X$.
\begin{define}
{\it A covering family of curves} on a projective variety $X$ is a diagram 
\[
\xymatrix{
 \Ca{C} \ar[r]^{f} \ar[d]_{\pi}& X\\
 T &       	\\
}	
\]
satisfying the following:
\begin{itemize}
\item[(i)]$\Ca{C}$ and $T$ are irreducible varieties,
\item[(ii)]$\pi:\Ca{C} \to T$ is a smooth family of projective curves,
\item[(iii)]$f$ is a dominant morphism
\item[(iv)]for a general fiber $\Ca{C}_{t}$ ($t \in T$), the restriction morphism $f_{t}:\Ca{C}_{t} \to X$
is birational onto its image.
\end{itemize}
\end{define}
The basic properties of the covering family of curves are summarized in \cite[Remark~1.5]{BDELU}.

For a projective curve $C$, the {\it gonality} is defined by 
\[
\Rm{gon}(C):=\Rm{min}\{\Rm{deg}(\varphi) \mid \textrm{rational covering}\;\varphi:C \dashrightarrow \Bb{P}^1\}.
\]
It is well-known that the function $T \ni t \to \Rm{gon}(\Ca{C}_{t}) \in \Z_{\geq 0}$ is a lower semicontinuous function
for a smooth family $\pi: \Ca{C} \to T$ of projective curves.
A covering family ($\pi: \Ca{C} \to T$, $f: \Ca{C} \to X$) is {\it $c$-gonal} if for a general point $t \in T$,
it holds $\Rm{gon}(\Ca{C}_{t})=c$.

Let us introduce two birational invariants covering gonality and minimal covering degree,
as analogous to the gonality. 
\begin{define}[Gonality type invariants]
Consider a $n$-dimensional projective variety $X$. 

\noindent
(1)\;\textit{Covering gonality}
\[
\Rm{cov.gon}(X):=\Rm{min}\{c\mid \exists c\text{-gonal covering family of }X\}.
\]

\noindent
(2)\:\textit{Minimal covering degree}
\[
\Rm{mcd}(X):=\Rm{min}\{\Rm{deg}(\varphi) \mid\text{rational covering}\;\varphi:X \dashrightarrow Y \text{ with } \Rm{deg}(\varphi) \geq 2\}.
\]
Note that, in the definition of the minimal covering degree,
the "min" is taken by all projective varieties $Y$ and 
it takes the minimal degree of all rational coverings between $X$ and $Y$.
\end{define}

We can state our main results.
We show the slope inequality for fibered 3-folds, 
in which the covering gonality and minimal covering degree of the fibers are bounded.
\begin{thm}[Theorem~\ref{mainthm}]
Let $f:X \to B$ be a relatively minimal fibered $3$-fold and $c \geq 4$ be a positive integer.
Assume a general fiber $F$ of f is of general type with $p_{g}(F) \geq 1$.

If $\Rm{cov.gon}(F)$ and $\Rm{mcd}(F)$ are at least $c$ for a general fiber $F$,
then 
\[
K_{f}^3 \geq \left( 1 - \frac{4}{c}\right)\frac{K_{F}^2}{K_{F}^2 + 24\left( 1 - \frac{4}{c}\right)}\left(10 \cdot\Rm{deg}f_{\ast}\omega_{f} -2l(2)\right).
\] 
Here, $l(2)$ is the correction term in the Reid-Fletcher's plurigenera formula $($see \cite[Definition~2.6]{Fl}$)$.
\end{thm}

The covering gonality has been studied by Bastianelli as a measure of rationality for projective varieties (\cite{Bas}). 
The values of the covering gonality for various projective varieties have been investigated. 
For instance, Bastianelli has studied the relationship between the gonality of a smooth projective curve 
and the covering gonality of its second symmetric product (\cite{Bas}). 
Bastianelli-De Poi-Ein-Lazarsfeld-Ullery have studied the relationship between the degree of hypersurfaces in projective spaces and their covering gonality (\cite{BDELU}).

The minimal covering degree is a new birational invariant introduced by us.
It is an analogy for the degree of irrationality.
For a $n$-dimensional projective varieties $X$, the \textit{degree of irrationality} is defined by 
\[
\Rm{irr}(X):=\Rm{min}\{\Rm{deg}(\varphi) \mid\text{rational covering}\;\varphi:X \dashrightarrow \Bb{P}^{n}\}.
\]

\subsection{Background of the geography of fibered spaces}
For a fiber space $f:X \to B$ let $\omega_{f}$ be the relative dualizing sheaf and $K_f$ relative canonical divisor.
It holds that $\Ca{O}_{X}(K_{f}) = \omega_{f}$ (\cite[Proposition~5.75]{KM}).
Then we consider the relative numerical invariants $K_{f}^{\Rm{dim}(X)}$, $\Rm{deg}(f_{\ast}\omega_{f})$
and $\chi_{f}:=(-1)^{\Rm{dim}(X)}(\chi(\Ca{O}_{X})-\chi(\Ca{O}_{F})\chi(\Ca{O}_{B}))$.
We note that $\Rm{deg}(f_{\ast}\omega_{f}) \geq 0$ by Fujita's nefness Theorem.
Furthermore, if $f:X \to B$ is relatively minimal,
that is, $X$ is a $\Bb{Q}$-factorial variety with only terminal singularities 
and the canonical divisor $K_{X}$ is $f$-nef,
then $K_{f}^{\Rm{dim}(X)} \geq 0$ (\cite[Theorem~1.4]{Ohno}).
If $\Rm{deg}f_{\ast}\omega_{f} > 0$,
then the {\it slope} of $f:X \to B$ is defined as the following positive ratio: 
\[
\lambda(f) := \frac{K_{f}^{\Rm{dim}(X)}}{\Rm{deg}(f_{\ast}\omega_{f})}
\]
Evaluating the value of the slope for fibered spaces 
is a fundamental problem in the geography of fibered spaces. 
The slope inequalities contribute to understanding the relationship between the geometric properties of fibers and relative numerical invariants.

Regarding fibered surfaces, various slope inequalities have been established.
Let $f:X \to B$ be a fibered surface with genus $g(F) \geq 2$.
Among others, one most fundamental results is the following sharp lower bound:
\begin{align}
\label{slope}
\lambda(f) \geq \frac{4(g-1)}{g}
\end{align}
Xiao showed (\ref{slope}) by using a Harder-Narasimhan filtration method (\cite{Xiao}).

When equality in the inequality (\ref{slope}) holds, 
it is known that the general fiber of $f$ becomes a hyperelliptic curve. 
Hyperelliptic curves are special curves among algebraic curves. 
Therefore, the establishment of the slope inequality for fiber surfaces with more "general" curves as fibers becomes a crucial problem.
Reid showed the slope inequality $\lambda(f) \geq 3$
for non-hyperelliptic fibrations with genus $3$.
The slope inequality (\ref{slope}) implies $\lambda(f) \geq 8/3$ for fibered surfaces with genus $3$.
Therefore, Reid's slope inequality is stronger than (\ref{slope}) for genus $3$ non-hyperelliptic fibrations.

We focus on Lu-Zuo's slope inequality.
Lu-Zuo (\cite{Lu}) showed the slope inequality for fibered surfaces 
of which general fibers are minimal covering degree $\geq c$.
Our main result Theorem~\ref{mainthm} corresponds to a $3$-dimensional version of Lu-Zuo's slope inequality Theorem~\ref{Lu-Zuo}.
\begin{thm}[Lu-Zuo \cite{Lu}]
\label{Lu-Zuo}
Let $f:X \to B$ be a non-locally trivial and relatively minimal fibered surface with genus $g(F) \geq 2$.
If $\Rm{mcd}(F) = c \geq 3$ for a general fiber $F$,
then 
\begin{align*}
\lambda(f) > \left(5-\frac{1}{c-1}\right)\left(\frac{g-1}{g+2}\right).
\end{align*}
\end{thm}

For fibered $3$-folds,
Ohno (\cite{Ohno}), Hu-Zhang (\cite{Hu}), and others have been exploring Xiao's type slope inequality (\ref{slope}). 
On the other hand, studies on slope inequalities for fibered 3-folds with specific geometric properties endowed to the fibers, such as the approach of Reid and Lu-Zuo, are limited.
Therefore, we have established a $3$-dimensional version of the slope inequality shown by Lu-Zuo.

From the perspective of the moduli theory of curves, we interpret Lu-Zuo's slope inequality (\ref{Lu-Zuo}). 
Let $M_{g}$ be the coarse moduli space of smooth projective curves of genus $g \geq 2$. 
Fix a Zariski closed subset $Z \subset M_{g}$ of codimension at least $1$. 
We refer to a fiber surface $f: S \to B$ with a general fiber being in $M_{g} \setminus Z$ as a \textit{$Z$-general fibered surface} (cf.\cite{Eno2}). 
One can pose the following question Problem~\ref{prob2} related to the slope inequality.
This problem is related to the Farkas-Morrison conjecture (\cite[Conjecture~0.1]{HM}) for the moduli space.

\begin{prob}
\label{prob2}
Let $f: X \to B$ be a $Z$-general fibered surface.
Assume $f$ is non-locally trivial and relatively minimal.
Establish the slope inequality for $f$:
\[ \lambda(f) \geq \lambda_{(Z)} \]
where $ \lambda_{(Z)} $ is a positive rational number determined by the closed subset $Z \subset M_{g}$.
\end{prob}

Fix an integer $c \geq 3$. The subset of the moduli space
\[ \mathcal{C}_{g}(c) := \{x \in M_{g} \mid \mathrm{mcd}(C_{x}) \leq c-1 \} \]
is a Zariski closed subset (cf. \cite[Theorem~8.23]{ACG}). 
Similar to gonality, 
the minimal covering degree also exhibits "well" behavior on the moduli space of curves, 
as this subset is Zariski-closed. 
Furthermore, when $c-1 < \lfloor \frac{g+3}{2} \rfloor$, 
we have $ \mathrm{dim}\;\mathcal{C}_{g}(c) < \mathrm{dim}\;M_{g} $. 
Hence, Lu-Zuo's slope inequality can be interpreted as follows:

\begin{thm}[=Theorem~\ref{Lu-Zuo}]
Let $f:X \to B$ be a non-locally trivial and relatively minimal fibered surface with genus $g(F) \geq 2$.
If $f$ is $ \mathcal{C}_{g}(c)$-general for $c \geq 3$, then the slope inequality
\[ \lambda(f) > \left(5-\frac{1}{c-1}\right)\left(\frac{g-1}{g+2}\right) \]
holds.
\end{thm}

Our study of the slope inequality for fibered 3-folds in this paper is motivated by the perspective of incorporating Problem~\ref{prob2} into the geographical study of fibered $3$-folds.
To incorporate the perspective of Problem~\ref{prob2} into the geographical study of fibered $3$-folds, 
it is essential to investigate whether the covering gonality and minimal covering degree behave "well" on the moduli space of surfaces. 
However, unlike in the case of projective curves, 
this behavior is still unknown. 
For example, it is unclear whether the covering gonality and minimal covering degree exhibit lower semi-continuity, and so on. 
Therefore, unfortunately, as of now, 
the moduli-theoretic interpretation of Theorem~\ref{mainthm} has not been provided.
This is a subject for future study.

\vspace{3\baselineskip}

\textbf{Acknowledgememt}

The author is deeply grateful to Professor Kazuhiro Konno and Professor Tadashi Ashikaga 
for his valuable advice and support.
The author also thanks Doctor Makoto Enokizono for his precious advice.

\setcounter{tocdepth}{1}
\tableofcontents

\section{Preliminaries}
We summarize the basic tools used in the proof of the main result in this section.
\subsection{Nef threshold of $\Bb{Q}$-twisted locally free sheave}
For the convenience of the reader,
we recall the basic properties of $\Bb{Q}$-twisted locally free sheaves.

Let $B$ be a smooth projective curve.
For a locally free sheaf $\Ca{E}$ and a $\Bb{Q}$-divisor $\delta$,
we denote by $\Ca{E}\langle\delta\rangle$ the $\Bb{Q}$-twist of $\Ca{E}$ by $\delta$.
We recall the basic notions of $\Bb{Q}$-twist $\Ca{E}\langle\delta\rangle$:
\begin{itemize}
\item Define $\Rm{deg}(\Ca{E}\langle\delta\rangle):=\Rm{deg}(\Ca{E})+\Rm{rk}(\Ca{E})\cdot \Rm{deg}(\delta)$.
\item Let $\pi : \Bb{P}(\Ca{E}) \to B$ be the projective bundle associated to $\Ca{E}$.
A $\Bb{Q}$-twisted locally free sheaf $\Ca{E}\langle\delta\rangle$ is nef if the 
$\Bb{Q}$-invertible sheaf $\Ca{O}_{\Bb{P}(\Ca{E})}(1)\otimes \pi^{\ast} \delta $ is nef.
\item Put $\mu_{\Rm{nef}}(\Ca{E}):=\Rm{max} \{\Rm{deg}(\delta)\mid \Ca{E}\langle -\delta \rangle \textrm{ is nef }\}$.
\end{itemize}
\begin{define}
Let $\Ca{E}$ be a locally free sheaf on $B$.
Considering the Harder-Narasimhan filtration $0=\Ca{E}_0 \subset \Ca{E}_1 \subset \cdots \subset  \Ca{E}_{n-1} \subset \Ca{E}_n = \Ca{E}$,
put $\mu_{\Rm{min}}(\Ca{E}):=\mu(\Ca{E}_{n}/\Ca{E}_{n-1})$.
\end{define}

\begin{lem}
\label{nefth}
If $\Ca{E}$ is a locally free sheaf on $B$,
then it holds that $\mu_{\Rm{nef}} (\Ca{E}) = \mu_{\Rm{min}} (\Ca{E})$.
\end{lem}
\begin{proof}
Let $\delta$ be a $\Bb{Q}$-divisor on $B$.
By \cite[Corollary~3.8]{Nak}, $\Bb{Q}$-line bundle $\Ca{O}_{\Bb{P}(\Ca{E})}(1)\otimes \pi^{\ast} \delta $ is nef 
if and only if  $\Rm{deg}(\delta) \leq  \mu(\Ca{E}_{n}/\Ca{E}_{n-1})$.
Therefore, we have $\mu_{\Rm{nef}} (\Ca{E}) = \mu(\Ca{E}_{n}/\Ca{E}_{n-1})$.
\end{proof}

One can show the following lemma. The proof is straightforward and will be omitted.
\begin{lem}
\label{easy}
 Let $\Ca{F} \to \Ca{G}$ be a surjective map between locally free sheaves on $B$ and 
 $\delta$ be a $\Bb{Q}$-divisor on $B$.
 If $\Ca{F}\langle\delta\rangle$ is nef, then $\Ca{G}\langle\delta\rangle$ is also nef.
\end{lem}

\begin{lem}
\label{tensor}
Let $\Ca{E}\langle\delta\rangle$ be a  $\Bb{Q}$-twist of $\Ca{E}$ by $\delta$.
Then the following are equivalent.
\begin{itemize}
\item[(1)] There exists a positive integer $m \in \Z$ such that $\Ca{E}^{\otimes m} \langle m\delta \rangle$ is nef.
\item[(2)] The $\Bb{Q}$-twist $\Ca{E}\langle\delta\rangle$ is nef.
\end{itemize}
\end{lem}
\begin{proof}
$(1)\Rightarrow(2)$: Let $m$ be a positive integer $m \in \Z$ such that $\Ca{E}^{\otimes m} \langle m\delta \rangle$ is nef.
Since there exists the natural surjective map $\Ca{E}^{\otimes m} \twoheadrightarrow \Rm{Sym}^{m}(\Ca{E})$ 
and $\Ca{E}^{\otimes m} \langle m\delta \rangle$ is nef,
$\Bb{Q}$-twist $(\Rm{Sym}^{m}(\Ca{E}))\langle m\delta \rangle$ is also nef by Lemma~\ref{easy}.
It means that $\Ca{O}_{\Bb{P}(\Rm{Sym}^{m}(\Ca{E}))}(1)\otimes \pi^{\ast}_{m} (m\delta) $ is a nef invertible sheaf on $\pi_{m} : \Bb{P}(\Rm{Sym}^{m}(\Ca{E})) \to B$.
Let $\pi : \Bb{P}(\Ca{E}) \to B$ be the projective bundle associated to $\Ca{E}$.
Consider the relative $m$-th Veronese embedding:
\[
\xymatrix{
 \Bb{P}(\Ca{E}) \ar@{^{(}-_>}[rr] \ar[rd]_{\pi}& &\Bb{P}(\Rm{Sym}^{m}(\Ca{E})) \ar[dl]^{\pi_{m}}\\
 & B &     	\\
}	
\]
The pullback of $\Ca{O}_{\Bb{P}(\Rm{Sym}^{m}(\Ca{E}))}(1)\otimes \pi^{\ast}_{m} (m\delta) $ by the Veronese embedding
is $\Ca{O}_{\Bb{P}(\Ca{E})}(m)\otimes \pi^{\ast}( m\delta )$.
Hence $\Ca{O}_{\Bb{P}(\Ca{E})}(m)\otimes \pi^{\ast} (m\delta) $ is nef,
and it implies that $\Ca{O}_{\Bb{P}(\Ca{E})}(1)\otimes \pi^{\ast} \delta $ is nef.

$(2)\Rightarrow(1)$: See \cite[Theorem~6.2.12 (iii)]{Laz}.
\end{proof}

\begin{lem}
\label{nefineq}
The followings hold.
\begin{itemize}
\item[(1)] Let $\Ca{F} \to \Ca{G}$ be a surjective map between locally free sheaves on $B$.
Then, 
\[
\mu_{\Rm{nef}}(\Ca{G}) \geq \mu_{\Rm{nef}}(\Ca{F}).
\]
\item[(2)] Let $\Ca{E}$ be a locally free sheaf on $B$. Then, 
\[
\mu_{\Rm{nef}}(\Ca{E}\otimes\Ca{E}) = 2 \mu_{\Rm{nef}}(\Ca{E}).
\]
\end{itemize}
\end{lem}
\begin{proof}
$(1)$: Let $\pi_{f} : \Bb{P}(\Ca{F}) \to B$ and $\pi_{g} : \Bb{P}(\Ca{G}) \to B$ be projective bundles associated to $\Ca{F}$ and $\Ca{G}$, repectively.
We denote by $\Gamma_f$ and $\Gamma_g$ a fiber of $\pi_f$ and $\pi_g$, respectively.
Then it suffices to show 
\[
\Rm{max} \{t \mid H_f - t \Gamma_f  \textrm{ is nef }\} 
\leq
\Rm{max} \{t \mid H_g - t \Gamma_g  \textrm{ is nef }\},
\]
where $H_f \in |\Ca{O}_{\Bb{P}(\Ca{F})}(1)|$ and $H_g \in |\Ca{O}_{\Bb{P}(\Ca{G})}(1)|$.
Now, there exists a morphism $j:\Bb{P}(\Ca{G}) \to \Bb{P}(\Ca{F})$ over $B$ corresponding to the surjective map $\Ca{F} \to \Ca{G}$.
Considering the pullback by $j$, we have 
\[
\{t \mid H_f - t \Gamma_f  \textrm{ is nef }\} 
\subset
\{t \mid H_g - t \Gamma_g  \textrm{ is nef }\}.
\]
Therefore, we get the desired inequality.

$(2)$: We can show as follows.
\begin{align*}
\mu_{\Rm{nef}}(\Ca{E}\otimes\Ca{E}) 
&= \Rm{max} \{\Rm{deg}(\delta)\mid (\Ca{E}\otimes \Ca{E})\langle -\delta \rangle \textrm{ is nef }\} \\
&=\Rm{max} \{2\cdot\Rm{deg}(\delta)\mid (\Ca{E}\otimes \Ca{E})\langle -2\delta \rangle \textrm{ is nef }\} \\
&= \Rm{max} \{2\cdot\Rm{deg}(\delta)\mid \Ca{E}\langle -\delta \rangle \textrm{ is nef }\} \quad  (\textrm{by Lemma~} \ref{tensor},(2))\\
&= 2 \mu_{\Rm{nef}}(\Ca{E})
\end{align*}
\end{proof}

\subsection{Clifford Plus Theorem}
We state the Clifford Plus theorem for projective surfaces in this subsection. 
First, we introduce the classical Clifford Plus Theorem for projective curves.
\begin{prop}[Clifford Plus Theorem]
\label{CPthm1}
Let $C$ be a projective curve and let $\Ca{L}$ be an invertible sheaf on $C$.
Consider a $\C$-subspace $V \subset H^{0}(C,\Ca{L})$ with $\Rm{dim}_{\C}V=r+1$.
Assume that the corresponding morphism $\varphi_{V}:C \to \Bb{P}(V)$ is a closed immersion.
Then 
\[
\Rm{rank}\left(V\otimes V \to  H^{0}(C,\Ca{L}^2) \right) \geq \Rm{min}\{3r, p_{a}(C)+2r+1\},
\]
where $V\otimes V \to  H^{0}(C,\Ca{L}^2)$ is the multiplicative map.
\end{prop}
The proof of Theorem~\ref{CPthm1} can be found in \cite[\S III.2]{ACGH}.
We refer to the following inequality that appears in the proof of the Clifford Plus Theorem.

\begin{rmk}
\label{preCP}
Let the notations and assumptions be as above.
Then, it holds that
\[
\Rm{rank}\left(V\otimes V \to  H^{0}(C,\Ca{L}^2) \right) \geq \Rm{min}\{3r, \Rm{deg}(C)+r+1\},
\]
\end{rmk}

Theorem\ref{CPthm1} is extended to the case of surfaces as the following.
The essential part of it has been shown by Konno in \cite[Lemma~1.2]{Ko}.

\begin{prop}
\label{CPthm2}
Let $S$ be a projective surface (not necessarily normal) and let $\Ca{L}$ be an invertible sheaf on $S$.
Denote by $p_{g}(S)$ a geometric genus of a smooth model of $S$.
Consider a $\C$-subspace $V \subset H^{0}(S,\Ca{L})$ with $\Rm{dim}_{\C}V=r+1$.
Assume that the corresponding morphism $\varphi_{V}:S \to \Bb{P}(V)$ is a closed immersion.
Then 
\begin{equation}
\nonumber
\Rm{rank}\left(V\otimes V \to  H^{0}(S,\Ca{L}^2) \right) \geq 
 \begin{cases}
   4r-2 & \text{if $p_{g}(S) \geq 1$,} \\
      3r       & \text{if $p_{g}(S)=0$,}
  \end{cases}
\end{equation}
where $V\otimes V \to  H^{0}(S,\Ca{L}^2)$ is the multiplicative map.
\end{prop}
\begin{proof}
By \cite[Lemma~1.2]{Ko}, we have 
\begin{equation}
\label{9/14,2}
\Rm{rank}\left(V\otimes V \to  H^{0}(S,\Ca{L}^2) \right) \geq \Rm{min}\{4r-2, \Rm{deg}(S)+2r+1\}.
\end{equation}

We will show that
\begin{equation}
\label{9/14,1}
\Rm{deg}(S) \geq
\begin{cases}
   \Rm{min}\{3r-5,p_{g}(S)+2r-3 \} & \text{if $p_{g}(S) \geq 1$,} \\
      r-1     & \text{if $p_{g}(S)=0$.}
  \end{cases}
\end{equation}
Since $S$ is non-degenerate,
it holds $\Rm{deg}(S) \geq r-1$.
Hence it suffices to show $\Rm{deg}(S) \geq \Rm{min}\{3r-5,p_{g}(S)+2r-3 \}$
if $p_{g}(S)\geq 1$.
Assume $3r-5\geq \Rm{deg}(S)$.
Let $m = \left[\frac{\Rm{deg}(S)-1}{r-2} \right]$ and let $\varepsilon = (\Rm{deg}(S)-1) - m(r-2) \geq 0$,
then $m=0,1,2$ or $3$ by $\varepsilon \geq 0$.
By Castelnuovo-Harris bound \cite[\S 2]{Harris}, we have 
\begin{align}
\label{deg2(1)}
p_{g}(S) \leq \binom{m}{3}(r-2)+\frac{1}{2}m(m-1)\varepsilon,
\end{align}
where the binomial coefficient $\binom{a}{b}=0$ when $b>a$.
Since the assumption $p_{g}(S) \geq 1$ and (\ref{deg2(1)}),
we have $m=2$ or $3$.
If $m=2$, then $p_{g}(S)\leq \Rm{deg}(S) -2r + 3$ by (\ref{deg2(1)}).
If $m=3$, then $\varepsilon = (\Rm{deg}(S)-1) -3r +6 \geq 0$.
By the assumption $\Rm{deg}(S)\leq 3r-5$,
we get $\Rm{deg}(S)=3r-5$.
Hence we have (\ref{9/14,1}).
By (\ref{9/14,2}) and (\ref{9/14,1}), we get the desired inequality.
\end{proof}



\subsection{The properties of the covering gonality and minimal covering degree}

We introduce the property of the covering gonality used in the proof of the main result.
\begin{prop}
\label{gonality}
Let $S$ be a minimal surface of general type with $p_{g}(S) \geq 2$.
Let $V \subset H^{0}(K_{S})$ be a subspace with $\Rm{dim}V=r+1$, 
let $\varphi_{V}: S \dashrightarrow \Bb{P}V$ be the rational map corresponding to $V$,
let $\Gamma$ be a Zariski closure of image of $\varphi_{V}$ and
let $\mu: \Ti{S} \to S$ be a resolution of the indeterminacy locus of $\varphi_{V}$.
Denote by $\iota_{V}:\Ti{S} \to \Gamma$, $\Ti{\iota}_{V}:\Ti{S} \to \Ti{\Gamma}$ the induced morphism,
the stein factorization of $\iota_{V}$, respectively.
\[
\xymatrix{
\Ti{S}  \ar[d]_{\mu}  \ar[rd]^{\iota_{V}} &  \\
S    \ar@{-->}[r]_{\varphi_{V}} &    \Gamma. \\
}	
\]
Take $M \in |\iota_{V}^{\ast}\Ca{O}_{\Gamma}(1)|$.
If $\Rm{dim}(\Gamma) = 1$ and $\Rm{cov.gon}(S)\geq c$ for some positive integer $c \geq 2$,
then
\[
(\mu^{\ast}K_{S}\cdot M) \geq 2r(c-2).
\]
\end{prop}
\begin{proof}
The morphisms  $\iota_{V}:\Ti{S} \to \Gamma$ and $\Ti{\iota}_{V}:\Ti{S} \to \Ti{\Gamma}$
fit the following diagram:
\[
\xymatrix{
  & \Ti{\Gamma} \ar[d]  \\
\Ti{S}   \ar[r]_{\iota_{V}} \ar[ru]^{\Ti{\iota}_{V}} &    \Gamma \\
}	
\]
In the above diagram, let $\beta$ be the degree of morphism $\Ti{\Gamma} \to \Gamma$.
Taking a general hyperplane $H \in |\Ca{O}_{\Bb{P}V}(1)|$,
we may assume $H|_{\Gamma}$ is reduced divisor
\[
 H|_{\Gamma}= p_1 + p_2 +\cdots +p_{\alpha}
\]
where $\alpha = \Rm{deg}\Gamma$ is a positive integer.
Since $\Gamma$ is non-degenerate, we have $\alpha \geq r$.
Let 
\[
(\iota_{V})^{\ast} H|_{\Gamma} = \sum_{\substack{i=1,\cdots,\beta \\ j=1, \cdots, \alpha}}\Ti{G}_{i,j} \in |\iota_{V}^{\ast}\Ca{O}_{\Gamma}(1)|,
\]
where $\sum_{i=1}^{\beta}\Ti{G}_{i,j}$ is a fiber of $\iota_{V}$ over $p_j$.
If necessary, we retake $H \in |\Ca{O}_{\Bb{P}V}(1)|$,
we may assume each $\Ti{G}_{i,j}$ is a smooth irreducible curve.

Since $\mu^{\ast}K_{S} \geq M$, there exists an effective divisor $E \in |\mu^{\ast}K_{S}- \sum \Ti{G}_{i,j}|$. 
So we have $\mu_{\ast}E \in |K_{S}- \sum G_{i,j}|$, where $G_{i,j}:=\mu_{\ast}\Ti{G}_{i,j}$.
Furthermore, since $G_{i,j}^2 \geq 0$ for any $i,j$,
we have $(K_{S}- \sum_{i,j} G_{i,j}) \sum_{i,j} G_{i,j} \geq 0$.
Therefore, we get 
\begin{align*}
2K_{S}\cdot \sum_{i,j}G_{i,j} & \geq K_{S}\cdot \sum_{i,j} G_{i,j}+ \left(\sum G_{i,j}\right)^2\\
& = \left(K_{S}+\sum_{i,j} G_{i,j}\right) \cdot \sum_{i,j} G_{i,j}\\
& = 2p_{a}\left(\sum_{i,j} G_{i,j}\right)-2\\
& \geq  \sum_{i,j} \left( 2p_{a}\left(G_{i,j}\right)-2 \right) \\
& \geq  \sum_{i,j}\left( 2g\left(\Ti{G}_{i,j}\right)-2 \right) \\
& \geq  r \left( 2g\left(\Ti{G}_{i,j}\right)-2 \right) \quad (\textrm{by } \alpha \geq r). 
\end{align*}
We have 
\[
\left[\frac{g\left(\Ti{G}_{i,j}\right) +3}{2} \right] \geq \Rm{gon}\left(\Ti{G}_{i,j}\right) \]
by \cite[\S V, Theorem~1.5]{ACGH}.
Since the gonality of general fibers of $\Ti{\iota}_{V}$ is at least a positive integer $c\geq 2$, we have
\[
\Rm{gon}\left(\Ti{G}_{i,j}\right)  \geq c.
\]
Therefore, we have 
\[
2K_{S}\cdot \sum_{i,j}G_{i,j} \geq 4r(c-2).
\]
On the other hand, we have 
$(K_{S}\cdot \sum_{i,j}G_{i,j}) = (\mu^{\ast}K_{S}\cdot \sum_{i,j}\Ti{G}_{i,j}) = (\mu^{\ast}K_{S}\cdot M )$.
Hence we get 
\[
(\mu^{\ast}K_{S}\cdot M) \geq 2r(c-2).
\]
\end{proof}

\begin{lem}
\label{gcg}
Let $S$ be a projective surface. 
Take an invertible sheaf $\Ca{L} \in \Rm{Pic}(S)$.
and consider a subspace $V \subset H^0(S,\Ca{L})$ with $\Rm{dim}V = r+1$.
\[
\xymatrix{
\Ti{S}  \ar[d]_{\mu}  \ar[rd]^{\iota_{V}} & & \\
S    \ar@{-->}[r]_{\varphi_{V}} &    \Gamma  \ar@{^{(}-_>}[r] & \Bb{P}(V)\\
}	
\]
Assume $\iota_{V}$ is a generically finite morphism with degree $\geq 2$.
If $\Rm{mcd}(S) \geq c$ for some positive integer $c \geq 2$,
then $\Rm{deg}(\iota_{V}) \geq c$.
\end{lem}
\begin{proof}
It holds clearly by the definition of the minimal covering degree.
\end{proof}

\begin{exa}[Surfaces with large covering gonality and minimal covering degree]
Fix an arbitrary positive integer $m \in \Z_{\geq 1}$.
We construct surfaces with covering gonality and minimal covering degree $\geq m$.
Denote by $C^{(2)}$ the $2$-th symmetric product of a smooth projective curve $C$.

To construct such surfaces, we introduce results by Bastianelli and Lee-Pirola.
Lee-Pirola (\cite[Theorem~4.6]{LP}) show that 
if $C$ is a very general curve of genus $g \geq 10$,
then 
\[
\Rm{irr}(C^{(2)})=\Rm{mcd}(C^{(2)}).
\]
Bastianelli (\cite[Theorem~1.4, 1.6]{Bas}) show the following:
Let $C$ be a smooth projective curve of genus $g \geq 7$.
\begin{itemize}
\item[(1)] It holds that $\Rm{cov.gon}(C^{(2)}) = \Rm{gon}(C)$.
\item[(2)] If $\Rm{gon}(C) \geq 6$, then it holds that $\Rm{irr}(C^{(2)}) \geq \Rm{gon}(C)$.
\end{itemize}
From the above two results, for a very general curve with $\Rm{gon}(C) \geq 6$,
it holds that 
\[
\mathrm{min}\{\Rm{irr}(C^{(2)}), \Rm{mcd}(C^{(2)})\} \geq \Rm{gon}(C).
\]
Since there exists a smooth projective curve $C$ with gonality $\geq m$,
there exists a smooth projective surface $C^{(2)}$ with covering gonality and minimal covering degree $\geq m$.
\end{exa}

\section{Slope of fibered $3$-folds}
In this section, we present the proof of the main result, Theorem~\ref{mainthm}. 
We consider a relatively minimal fibered $3$-folds $f:X \to B$ such that general fibers are of general type.
Assume geometric genus $p_g \geq 1$ of the general fiber.
The key elements for establishing Theorem~\ref{mainthm} are Propositions~\ref{key1} and \ref{key2}. 
In this section, we show these two propositions and use them to prove Theorem~\ref{mainthm}.
Proposition~\ref{key2} is where the assumptions regarding the covering gonality and minimal covering degree come into play.
First, we aim to prove Proposition~\ref{key1}. 
Subsequently, we will focus on the proof of Proposition~\ref{key2}.

We consider the Harder-Narasimhan filtration of $ f_{\ast}\omega_{f}$:
 \[
0=\Ca{E}_0 \subset \Ca{E}_1 \subset \Ca{E}_2 \subset \cdots \subset \Ca{E}_n =  f_{\ast}\omega_{f}.
\]
\begin{prop}[\cite{Ohno}, \S1]
\label{Ohno1}
There exist a smooth projective variety $Y$, a birational morphism $\mu: Y \to X$, and rational sections $\varphi_{i}:X \dashrightarrow \Bb{P}_{B}(\Ca{E}_{i})$ satisfying the following conditions.

\noindent
$(1)$ The morphism $\mu: Y \to X$ is a resolution of indeterminacy locus of
a rational section $\varphi_{i}$  induced by the natural map $f^{\ast}\Ca{E}_{i} \to \Ca{O}_{X}(K_{f})$ for all $i$. 
\begin{align}
\label{6/5,1}
\xymatrix{
Y  \ar[d]_{\mu}  \ar[rd]^{\lambda_{i}} &   \\
X    \ar@{-->}[r]_{\varphi_{i}} &    \Bb{P}_{B}(\Ca{E}_{i})\\
}	
\end{align}

\noindent
$(2)$ Put $M_{i}:=\lambda_{i}^{\ast}\Ca{O}_{\Bb{P}(\Ca{E}_{i})}(1) \in \Rm{Div}(Y)$ for $i = 1, \cdots, n$.
Then there exist an effective divisor $Z_{i}$ on $X$ and $\mu$-exceptional divisor $E_i$ on $Y$ such that
\[
M_{i} \sim_{\Bb{Q}} \mu^{\ast}(K_{f}-Z_{i})-E_{i}.
\]
Furthermore, it holds 
\[
Z_1 \geq Z_2 \geq \cdots \geq Z_{n} \geq Z_{n+1}:=0.
\]
\end{prop}
For a closed point $p \in B$, let
\[
\xymatrix{
\Ti{F}_{p}  \ar[d]_{\mu|_{\Ti{F}_{p}}}  \ar[rd]^{\iota_{i}} &   \\
F    \ar@{-->}[r] &    \Bb{P}(V_{i}) \\
}	
\]
be the base change of (\ref{6/5,1}) to $p \in B$,
where $\Ti{F}_{p}= Y \times_{B} \Rm{Spec}\;\Bb{C}(p)$, $F_{p}= X \times_{B} \Rm{Spec}\;\Bb{C}(p)$
$V_{i}:=\Ca{E}_{i}\otimes \Bb{C}(p)$ and $\iota_{i}:=\lambda_{i}|_{\Ti{F}_{p}}$ for $i=1, \cdots n$.
Let $\Gamma_{i} \subset \Bb{P}(V_{i})$ be the image of $\Ti{F}_{p}$ by $\iota_{i}$. 
Then we get the diagram
\begin{align}
\label{6/5,2}
\xymatrix{
\Ti{F}_{p}  \ar[d]_{\mu}  \ar[rd]^{\iota_{i}} &  & \\
F    \ar@{-->}[r] &    \Gamma_{i}  \ar@{^{(}-_>}[r]_{\tau_{i}} & \Bb{P}(V_{i}). \\
}	
\end{align}
Put 
\begin{itemize}
\item[] $l_{2}:=\Rm{min} \{i \mid \iota_{i} \textrm{ is a birational map }\}$ and
\item[] $l_{1}:=\Rm{min} \{i \mid \iota_{i} \textrm{ is a generically finite morphism with degree}\geq 2\}$.
\end{itemize}
If there exists no $i$ such that $\iota_{i}$ is a birational map, put $l_2 := n+1$.
If there exists no $i$ such that $\iota_{i}$ is a generically finite morphism with degree $\geq 2$, put $l_1 = l_2$.
Then we note that 
the morphism $\iota_{i}$ is birational if $l_{2} \leq i \leq n$, 
the morphism $\iota_{i}$ is generically finite with degree $\geq 2$ if $l_{1} \leq i \leq l_{2}-1$.

\begin{prop}
\label{prop2.6}
There exists a subsheaf $\Ca{F}_{i} \subset f_{\ast}\omega_{f}^{[2]}$ for $i=1, \cdots , n$ satisfying as follows:
\begin{itemize}
\item For $l_{2} \leq i \leq n$ $($i.e. $\iota_{i}$ is a birational morphism$)$,
it holds $\mu_{nef}(\Ca{F}_{i}) \geq 2 \mu_{i}$ and $\Rm{rank}(\Ca{F}_{i}) \geq 4 r_{i} -6 $.
\item For $l_{1} \leq i \leq l_{2}-1$ $($i.e. $\iota_{i}$ is a generically finite morphism$)$,
it holds $\mu_{nef}(\Ca{F}_{i}) \geq 2 \mu_{i}$ and $\Rm{rank}(\Ca{F}_{i}) \geq 3 r_{i} -3 $.
\item For $1 \leq i \leq l_{1}-1$, 
it holds $\mu_{nef}(\Ca{F}_{i}) \geq 2 \mu_{i}$ and $\Rm{rank}(\Ca{F}_{i}) \geq 2 r_{i} -1 $.
\end{itemize}
\end{prop}
\begin{proof}
For each $\Ca{E}_{i} \subset f_{\ast}\omega_{f}$, we consider the composition of natural map 
\[
\rho_{i}:\Ca{E}_{i} \otimes \Ca{E}_{i} \to \Rm{Sym}^{2}f_{\ast}\omega_{f} \to f_{\ast}\omega_{f}^{\otimes 2} \to f_{\ast}\omega_{f}^{[2]}.
\]
Let $\Ca{F}_{i}:=\Rm{Im}(\rho_{i}) \subset f_{\ast}\omega_{f}^{[2]}$.
We will show that $\Ca{F}_{i}$ satisfies the above conditions.

\noindent
\underline{Step~1}\;
We will show $\mu_{nef}(\Ca{F}_{i}) \geq 2 \mu_{i}$ for $i=1, \cdots, n$.
By Lemma~\ref{nefineq},
we have $\mu_{nef}(\Ca{F}_{i}) \geq \mu_{nef}(\Ca{E}_{i} \otimes \Ca{E}_{i}) = 2  \mu_{nef}(\Ca{E}_{i})$.
By Lemma~\ref{nefth}, we have $\mu_{nef}(\Ca{E}_{i})=\mu_i$.
Hence we get $\mu_{nef}(\Ca{F}_{i}) \geq 2 \mu_{i}$.

\noindent
\underline{Step~2}\;
We will give the estimation of $\Rm{rank}(\Ca{F}_{i})$.
Take a general point $p \in B$.
There exists the diagram (\ref{6/5,2}):
\[
\xymatrix{
\Ti{F}_{p}  \ar[d]_{\mu}  \ar[rd]^{\iota_{i}} &  & \\
F    \ar@{-->}[r] &    \Gamma_{i}  \ar@{^{(}-_>}[r]_{\tau_{i}} & \Bb{P}V_{i} \\
}	
\]
Then we have 
$
\Rm{rank}\left(V_{i} \otimes V_{i} \to H^{0}(\Gamma_{i}, \Ca{O}_{\Gamma_{i}}(2)) \right) = \Rm{rank}(\Ca{F}_{i}),
$
where $\Ca{O}_{\Gamma_{i}}(2) = \Ca{O}_{\Bb{P}(V_{i})}(2)|_{\Gamma_{i}}$.
Hence it suffices to estimate 
$
\Rm{rank}\left(V_{i} \otimes V_{i} \to H^{0}(\Gamma_{i}, \Ca{O}_{\Gamma_{i}}(2)) \right).
$

\noindent
$(1)$\;The case for $l_{2} \leq i \leq n$ $($i.e. $\iota_{i}$ is a birational morphism$)$.

\noindent
In this case, we have $p_{g}(S) \geq 1$.
By Proposition~\ref{CPthm2}, we get 
\[
\Rm{rank}\left(V_{i} \otimes V_{i} \to H^{0}(\Gamma_{i}, \Ca{O}_{\Gamma_{i}}(2)) \right) \geq 4r_{i}-6
\]
for $l_2 \leq i \leq n$.

\noindent
$(2)$\;The case for $l_{1} \leq i \leq l_{2}-1$ $($i.e. $\iota_{i}$ is a generically finite morphism with degree $\geq 2$$)$.

\noindent
By Proposition~\ref{CPthm2}, we get 
\[
\Rm{rank}\left(V_{i} \otimes V_{i} \to H^{0}(\Gamma_{i}, \Ca{O}_{\Gamma_{i}}(2)) \right) \geq 3r_{i}-3
\]
for $l_1 \leq i \leq l_{2}-1$.

\noindent
$(3)$\;The case for $1 \leq i \leq l_{1}-1$. 

\noindent
If $\Rm{dim}(\Gamma_{i})=1$, then we get 
\[
\Rm{rank}\left(V_{i} \otimes V_{i} \to H^{0}(\Gamma_{i}, \Ca{O}_{\Gamma_{i}}(2)) \right) \geq 2r_{i}-1
\]
by Theorem~\ref{CPthm1}.
If $\Rm{dim}(\Gamma_{i})=0$, then $i=1$ and $r_{1}=1$. Therefore it holds.
\end{proof}

\begin{prop}
\label{prop2.5}
Let $\Ti{\mu}_{1}>\cdots \Ti{\mu}_{k} \geq 0$ (resp. $0<\Ti{r}_{1}<\cdots < \Ti{r}_{k}$) 
be a decreasing (resp. increasing ) sequence of positive rational (resp. integer) numbers.
Assume that there exists locally free filtration $\Ca{F}_{1} \subset \cdots \subset \Ca{F}_{k}$ 
of $f_{\ast}\omega_{f}^{[2]}$ such that 
$\mu_{\Rm{nef}}(\Ca{F}_{i}) \geq \Ti{\mu}_{i}$ and $\Rm{rk}(\Ca{F}_{i}) \geq \Ti{r}_{i}$ for $i= 1, \cdots , k$.
Then,
\[
\frac{1}{2}K_{f}^3+3\chi_{f}+ l(2) \geq \sum_{i=1}^{k}\Ti{r}_{i}(\Ti{\mu}_{i}-\Ti{\mu}_{i+1}).
\]
Here, $l(2)$ is the correction term in the Reid-Fletcher's plurigenera formula $($see \cite[Definition~2.6]{Fl}$)$.
\end{prop}
\begin{proof}
Consider the Harder-Narasimhan filtration of $f_{\ast}\omega_{f}^{[2]}$:
\[
0=\Ca{E}'_0 \subset \Ca{E}'_1 \subset \Ca{E}'_2 \subset \cdots \subset \Ca{E}'_m = f_{\ast}\omega_{f}^{[2]}
\]
Put $\mu_{i}':=\mu(\Ca{E}'_{i}/\Ca{E}'_{i-1})$, $\mu'_{m+1}:=0$ 
and $r'_{i}:=\Rm{rk}(\Ca{E}'_{i})$ for $i=1, \cdots , m$. 
By \cite[Lemma~2.8]{Ohno}, we have 
\[
\frac{1}{2}K_{f}^3+3\chi_{f}+ l(2) = \Rm{deg}f_{\ast}\omega_{f}^{[2]} 
\]
Hence it suffices to show 
\[
\Rm{deg}f_{\ast}\omega_{f}^{[2]} \geq \sum_{i=1}^{m}\Ti{r}_{i}(\Ti{\mu}_{i}-\Ti{\mu}_{i+1}).
\]
To show it, we consider the Harder-Narashimhan polygon.
If we regard $\Ca{E}'_{i}$ as a point of $(\Rm{rk}(\Ca{E}'_{i}),\Rm{deg}(\Ca{E}'_{i}))$ in the two-dimensional coordinate system and
we connect points $\Ca{E}'_{i}$ and $\Ca{E}'_{i+1}$,
then we get the Harder-Narasimhan polygon.
An arbitrary subsheaf of $f_{\ast}\omega_{f}^{[2]}$ is in the interior or boundary of the Harder-Narasimhan polygon (\cite[Theorem~2]{Shatz}).
For $i=1,\dots, k$,
we denote by $(\Ca{F}_{i}/\Ca{F}_{i+1})_{free}$ the locally free part of $\Ca{F}_{i}/\Ca{F}_{i+1}$
and put $\mu(\Ca{F}_{i}/\Ca{F}_{i+1}):=\mu((\Ca{F}_{i}/\Ca{F}_{i+1})_{free})$.
Then, we get 
\begin{align*}
\sum_{i=1}^{m}r'_{i}(\mu'_{i}-\mu'_{i+1})  \geq \Rm{deg}(\Ca{F}_{k})
& \geq \sum_{i=1}^{k}\Rm{rk}(\Ca{F}_{i})\left(\mu(\Ca{F}_{i}/\Ca{F}_{i-1})-\mu(\Ca{F}_{i+1}/\Ca{F}_{i})\right) &\\
& = \sum_{i=1}^{k} \mu(\Ca{F}_{i}/\Ca{F}_{i-1}) \left(\Rm{rk}(\Ca{F}_{i})-\Rm{rk}(\Ca{F}_{i-1})\right),&
\end{align*}
where $\mu(\Ca{F}_{k+1}/\Ca{F}_{k}):=0$.
Since the convexity of the Harder-Narasimhan polygon of the locally free sheaf $(\Ca{F}_{i}/\Ca{F}_{i-1})_{free}$,
we have $\mu(\Ca{F}_{i}/\Ca{F}_{i-1}) \geq \mu_{\Rm{min}}(\Ca{F}_{i}/\Ca{F}_{i-1})$.
By Lemma~\ref{nefth},
we have $\mu(\Ca{F}_{i}/\Ca{F}_{i-1}) \geq \mu_{\Rm{min}}((\Ca{F}_{i}/\Ca{F}_{i+1})_{free})=\mu_{\Rm{nef}}((\Ca{F}_{i}/\Ca{F}_{i+1})_{free})$. 
Combing this with $\mu_{\Rm{nef}}((\Ca{F}_{i}/\Ca{F}_{i+1})_{free}) \geq \mu_{\Rm{nef}}(\Ca{F}_{i})$,
we have $\mu(\Ca{F}_{i}/\Ca{F}_{i-1}) \geq \mu_{\Rm{nef}}(\Ca{F}_{i})$.
Therefore we get 
\begin{align*}
 \sum_{i=1}^{k} \mu(\Ca{F}_{i}/\Ca{F}_{i+1}) \left(\Rm{rk}(\Ca{F}_{i})-\Rm{rk}(\Ca{F}_{i-1})\right)
& \geq \sum_{i=1}^{k} \mu_{\Rm{nef}}(\Ca{F}_{i}) \left(\Rm{rk}(\Ca{F}_{i})-\Rm{rk}(\Ca{F}_{i-1})\right) &\\
& \geq \sum_{i=1}^{k} \Ti{\mu}_{i} \left(\Rm{rk}(\Ca{F}_{i})-\Rm{rk}(\Ca{F}_{i-1})\right) &(\textrm{by } \mu_{\Rm{nef}}(\Ca{F}_{i}) \geq \Ti{\mu}_{i}) \\
& = \sum_{i=1}^{k}\Rm{rk}(\Ca{F}_{i})\left(\Ti{\mu}_{i}-\Ti{\mu}_{i+1}\right)   &\\
& \geq \sum_{i=1}^{k}\Ti{r}_{i} \left(\Ti{\mu}_{i}-\Ti{\mu}_{i+1}\right)  &(\textrm{by }\Ti{\mu}_{i}-\Ti{\mu}_{i+1} \geq 0).
\end{align*}
\end{proof}

\begin{prop}
\label{key1}
Let the notation be as above.
Then 
\[
K_{f}^3+2l(2) \geq \sum_{i=1}^{l_{1}-1}(2r_{i}-4)(\mu_{i}-\mu_{i+1})
+ \sum_{i=l_1}^{l_{2}-1}(6r_{i}-12)(\mu_{i}-\mu_{i+1})
+ \sum_{i=l_2}^{n}(10r_{i}-24)(\mu_{i}-\mu_{i+1}).
\]
\end{prop}
\begin{proof}
From Proposition~\ref{prop2.6}, for each $i=1, \cdots, n$, there exists the locally free sheaf 
$\Ca{F}_{i}$ such that 
\begin{align*}
\mu_{nef}(\Ca{F}_{i}) \geq 2\mu_{i} \quad \textrm{ and } \quad 
\Rm{rk}(\Ca{F}_{i}) \geq 
  \begin{cases}
    4r_{i}-6 & \text{if $l_2 \leq i \leq n$,} \\
    3r_{i}-3                 & \text{if $l_{1}\leq i \leq l_{2}-1$,} \\
    2r_{i}-1       & \text{if $1 \leq i \leq l_{1}-1$.}
  \end{cases}
\end{align*}
We can apply Proposition~\ref{prop2.6} to 
decreasing sequence 
\[
\{ 2\mu_{1}, \cdots , 2\mu_{l_{1}-1}, 2\mu_{l_{1}}, \cdots , 2\mu_{l_{2}-1}, 2\mu_{l_{2}}, \cdots , 2\mu_{n}\}
\]
and increasing sequence  
\[
\{ 2r_{1}-1, \cdots, 2r_{l_{1}-1}-1, 3r_{l_1}-3, \cdots , 3r_{l_{2}-1}-3, 4r_{l_{2}}-6, \cdots, 4r_{n}-6  \}.
\]
Therefore, we get 
\begin{align*}
&\frac{1}{2}K_{f}^3+3\chi_{f}+ l(2) \\
&\geq \sum_{i=1}^{l_{1}-1}(4r_{i}-2)(\mu_{i}-\mu_{i+1})
+ \sum_{i=l_1}^{l_{2}-1}(6r_{i}-6)(\mu_{i}-\mu_{i+1})
+ \sum_{i=l_2}^{n}(8r_{i}-12)(\mu_{i}-\mu_{i+1}).
\end{align*}
From \cite[Lemma~1.1]{Barja1}, we have 
\[
 \sum_{i=1}^{k}r_{i}(\mu_{i}-\mu_{i+1}) = \Rm{deg} f_{\ast}\omega_{f} \geq \chi_{f}.
\]
Hence we get 
\[
K_{f}^3+2l(2) \geq \sum_{i=1}^{l_{1}-1}(2r_{i}-4)(\mu_{i}-\mu_{i+1})
+ \sum_{i=l_1}^{l_{2}-1}(6r_{i}-12)(\mu_{i}-\mu_{i+1})
+ \sum_{i=l_2}^{n}(10r_{i}-24)(\mu_{i}-\mu_{i+1}).
\]
\end{proof}


From here, we focus on giving the proof of the Proposition~\ref{key2}.
To show Proposition~\ref{key2}, we use the following proposition by Ohno (\cite[p\;651]{Ohno}).

\begin{prop}$($\cite[p\;651]{Ohno}$)$
\label{Ohno2}
Let the notation and the assumption be as Proposition~\ref{Ohno1}.
Then, for a general fiber $\Ti{F}$ of $Y \to B$,
\[
K_{f}^3 \geq \sum_{i=1}^{n}(\mu_{i}-\mu_{i+1})\left((\mu|_{\Ti{F}}^{\ast}K_{F}\cdot M_{i}|_{\Ti{F}})+ (\mu|_{\Ti{F}}^{\ast}K_{F} \cdot M_{i+1}|_{\Ti{F}})\right).
\]
\end{prop}

\begin{prop}
\label{key2}
Let the notation and the assumption be as above.
Take a positive integer $c \geq 2$. 
If $\Rm{mcd}(F) \geq c$ for a general fiber $F$ of $f$,
then 
\[
K_{f}^3 \geq \sum_{i=1}^{l_{1}-1}(3c-4)(r_{i}-1)(\mu_{i}-\mu_{i+1})
+ \sum_{i=l_1}^{l_{2}-1}c(r_{i}-2)(\mu_{i}-\mu_{i+1}).
\]
\end{prop}
\begin{proof}


Replacing a general fiber $\Ti{F}$ if necessary, we may assume $F$ is smooth.
Put $d_{i}:=(\mu|_{\Ti{F}}^{\ast}K_{F}\cdot M_{i}|_{\Ti{F}})$.
From Proposition~\ref{Ohno2}, it suffices to show 
\begin{equation}
\label{6/6,1}
  d_{i} \geq 
 \begin{cases}
   0 & \text{if $l_2 \leq i \leq n$,} \\
    c(r_{i}-2)                 & \text{if $l_{1}\leq i \leq l_{2}-1$,} \\
      2(c-2)(r_{i}-1)       & \text{if $1 \leq i \leq l_{1}-1$.}
  \end{cases}
\end{equation}
Since $\mu|_{\Ti{F}}^{\ast}K_{F} \geq M_{i}|_{\Ti{F}}$,
we note that $d_{i}\geq M_{i}|_{\Ti{F}}^2$.
We recall the diagram (\ref{6/5,2}):
\[
\xymatrix{
\Ti{F}_{p}  \ar[d]_{\mu}  \ar[rd]^{\iota_{i}} &  & \\
F    \ar@{-->}[r] &    \Gamma_{i}  \ar@{^{(}-_>}[r]_{\tau_{i}} & \Bb{P}V_{i} \\
}	
\]

\noindent
$(1)$\;The case for $l_{2} \leq i \leq n$.

\noindent
Since $K_{F}$ is nef, we have $d_{i}\geq 0$.

\noindent
$(2)$\;The case for $l_{1} \leq i \leq l_{2}-1$.

\noindent
In this case, the morphism  $\iota_{i}: \Ti{F} \to \Gamma_{i}$ is a generically finite morphism with $\Rm{deg}(\iota_{i}) \geq 2$.
By $d_{i}\geq M_{i}|_{\Ti{F}}^2$, it suffices to show 
\[
M_{i}|_{\Ti{F}}^2 \geq c(r_{i}-2).
\]
By the assumption $\Rm{gcg}(\Ti{F}) \geq c$ and Lemma~\ref{gcg},
we have $\Rm{deg}(\iota_{i}) \geq c$.
Therefore, we get 
\[
d_{i} \geq M_{i}|_{\Ti{F}}^2 =\Rm{deg}(\iota_{i}) \cdot \Rm{deg}\Gamma_{i} \geq c(r_{i}-2).
\]
\noindent
$(3)$\;The case for $1 \leq i \leq l_{1}-1$.

\noindent
If $\Rm{dim}\Gamma_{i} = 1$,
then let $\Ti{\iota}_{i}: \Ti{F} \to \Ti{\Gamma}_{i}$ be the stein factorization of $\iota_{i}$:
\[
\xymatrix{
 &   \Ti{\Gamma}  \ar[d]  \\
\Ti{F}  \ar[ru]^{\Ti{\iota}_{i}}   \ar[r]_{\iota_{i}} &    \Gamma_{i}.   \\
}	
\]
By the definition $\Rm{gcg}(\Ti{F}) \geq c$,
a gonality of a general fiber of $\Ti{\iota}_{i}$ is at least $c$.
From Proposition~\ref{gonality}, we get 
\[
d_{i} =(\mu|_{\Ti{F}}^{\ast}K_{F}\cdot M_{i}|_{\Ti{F}}) \geq 2(c-2)(r_{i}-1).
\] 
Therefore, we get the desired inequalities (\ref{6/6,1}).
If $\Rm{dim}(\Gamma_{i})=0$, then $i=1$ and $r_{1}=1$. Therefore it holds.
\end{proof}

From here, we give the proof of the Theorem~\ref{mainthm}.

\begin{thm}
\label{mainthm}
Let $f:X \to B$ be a relatively minimal fibered $3$-fold and $c \geq 2$ be a positive integer.
Assume a general fiber $F$ of f is of general type.
If $\Rm{cov.gon}(F)$ and $\Rm{mcd}(F)$ are at least $c$ for a general fiber $F$,
then 
\[
K_{f}^3 \geq \left( 1 - \frac{4}{c}\right)\frac{K_{F}^2}{K_{F}^2 + 24\left( 1 - \frac{4}{c}\right)}\left(10 \cdot\Rm{deg}f_{\ast}\omega_{f} -2l(2)\right).
\] 
Here, $l(2)$ is the correction term in the Reid-Fletcher's plurigenera formula $($see \cite[Definition~2.6]{Fl}$)$.\end{thm}
\begin{proof}
By Proposition~\ref{key1} and \ref{key2},
we have 
\begin{align*}
K_{f}^3 + \frac{2(c-4)}{c}l(2)
&\geq \sum_{i=1}^{l_{1}-1}\left(\frac{c-4}{c}(2r_{i}-4)+\frac{4}{c}(3c-4)(r_{i}-1)\right)(\mu_{i}-\mu_{i+1})\\
&+ \sum_{i=l_1}^{l_{2}-1}\left(\frac{c-4}{c}(6r_{i}-12)+4(r_{i}-2)\right)(\mu_{i}-\mu_{i+1})\\
&+ \sum_{i=l_2}^{n}\frac{c-4}{c}\left(10r_{i}-24\right)(\mu_{i}-\mu_{i+1})\\
\end{align*}
For each $i$, one can show 
\begin{align*}
\left(1-\frac{4}{c}\right)(10r_{i}-24) \leq 
  \begin{cases}
   \left(1-\frac{4}{c}\right) \left(10r_{i}-24\right) & \text{if $l_2 \leq i \leq n$,} \\
  \left(1-\frac{4}{c}\right)(6r_{i}-12)+4(r_{i}-2)                & \text{if $l_{1}\leq i \leq l_{2}-1$,} \\
    \left(1-\frac{4}{c}\right)(2r_{i}-4)+\frac{4}{c}(3c-4)(r_{i}-1)       & \text{if $1 \leq i \leq l_{1}-1$.}
  \end{cases}
\end{align*}
Therefore we get 
\begin{align}
\label{main1}
K_{f}^3 + \frac{2(c-4)}{c}l(2) \geq 10\left(1-\frac{4}{c}\right)\Rm{deg}f_{\ast}\omega_{f}  -24 \left(1 - \frac{4}{c} \right)\mu_{1}
\end{align}

On the other hand, we have 
$K_{f}^3 \geq K_{F}^2 \mu_{1}$.
Indeed, a divisor $\mu^{\ast}K_{f}-Z_{1}-\mu_{1}\Ti{F}$ is nef by \cite[Lemma~2.3]{Ohno},
hence we have $(\mu^{\ast}K_{f}-Z_{1}-\mu_{1}\Ti{F})(\mu^{\ast}K_{f})^2 \geq 0$.
Therefore we have $K_{f}^3 \geq K_{F}^2 \mu_{1}$.
By (\ref{main1}) and $K_{f}^3 \geq K_{F}^2 \mu_{1}$, we get 
\[
K_{f}^3 \geq \left( 1 - \frac{4}{c}\right)\frac{K_{F}^2}{K_{F}^2 + 24\left( 1 - \frac{4}{c}\right)}\left(10 \cdot\Rm{deg}f_{\ast}\omega_{f} -2l(2)\right).
\] 
\end{proof}

\vspace{2\baselineskip}
{}
\bigskip
\bigskip

Hiroto Akaike

Department of Mathematics, 
Tohoku University,

e-mail: u802629d[at]alumni.osaka-u.ac.jp
\end{document}